\documentclass[a4paper,11pt,reqno]{amsart}

\usepackage{amssymb}

\newtheorem{theorem}[equation]{Theorem}
\newtheorem{lemma}[equation]{Lemma}
\newtheorem{corollary}[equation]{Corollary}

\numberwithin{equation}{section}

\theoremstyle{definition}
\newtheorem{definition}[equation]{Definition}
\newtheorem*{example*}{Example}

\newtheorem{remark}[equation]{Remark}
\newtheorem*{remark*}{Remark}

\newcommand{\bZ}{{\mathbb Z}}

\newcommand{\calV}{{\mathcal V}}
\newcommand{\frg}{{\mathfrak g}}
\newcommand{\frtg}{{\tilde{\mathfrak g}}}

\newcommand{\frf}{{\mathfrak f}}
\newcommand{\fre}{{\mathfrak e}}

\newcommand{\frs}{{\mathfrak s}}
\newcommand{\fra}{{\mathfrak a}}

\newcommand{\frgo} {{\frg_{\bar 0}}}
\newcommand{\frguno} {{\frg_{\bar 1}}}

\newcommand{\calS}{{\mathcal S}}

\newcommand{\calT}{{\mathcal T}}

\newcommand{\subo}{_{\bar 0}}
\newcommand{\subuno}{_{\bar 1}}

\newcommand{\bil}{{\textrm{b}}}

\providecommand{\espan}[1]{\text{span}\left\{ #1\right\}}

 \DeclareMathOperator{\tri}{\mathfrak{tri}}

 \DeclareMathOperator{\frosp}{\mathfrak{osp}}
 \DeclareMathOperator{\frsl}{{\mathfrak{sl}}}
  
 \DeclareMathOperator{\frsp}{{\mathfrak{sp}}}
 \DeclareMathOperator{\frso}{{\mathfrak{so}}}
 \DeclareMathOperator{\frpsl}{{\mathfrak{psl}}}
 \DeclareMathOperator{\frgl}{{\mathfrak{gl}}}
 \DeclareMathOperator{\frpgl}{{\mathfrak{pgl}}}

 \DeclareMathOperator{\ad}{ad}
 \DeclareMathOperator{\der}{\mathfrak{der}}
 \DeclareMathOperator{\inder}{\mathfrak{inder}}
 \DeclareMathOperator{\End}{End}
 
 \DeclareMathOperator{\Mat}{Mat}

\def\bigstrut{\vrule height 12pt width 0ptdepth 2pt}

\def\hregleta{\hrule height .5pt}
\def\hreglon{\hrule height1pt}
\def\vreglon{\vrule height 12pt width1pt depth 4pt}

\def\hreglonfill{\leaders\hreglon\hfill}

\def\hregletafill{\leaders\hregleta\hfill}

\newenvironment{romanenumerate}
 {\begin{enumerate}
 
 }{\end{enumerate}}

\begin{document}

\title{The extended Freudenthal Magic Square and Jordan algebras}

\author[Isabel Cunha]{Isabel Cunha$^{\diamond}$}
 \thanks{$^{\diamond}$ Supported by Centro de Matem\'atica da Universidade de Coimbra -- FCT}
 \address{Departamento de Matem\'atica, Universidade da Beira
 Interior,\newline 6200 Covilh\'a, Portugal}
 \email{icunha@mat.ubi.pt}

\author[Alberto Elduque]{Alberto Elduque$^{\star}$}
 \thanks{$^{\star}$ Supported by the Spanish Ministerio de
 Educaci\'on y Ciencia
 and FEDER (MTM 2004-081159-C04-02) and by the
Diputaci\'on General de Arag\'on (Grupo de Investigaci\'on de
\'Algebra)}
 \address{Departamento de Matem\'aticas, Universidad de
Zaragoza, 50009 Zaragoza, Spain}
 \email{elduque@unizar.es}

\date{\today}



\begin{abstract}
The Lie superalgebras in the extended Freudenthal Magic Square in
characteristic $3$ are shown to be related to some known simple Lie
superalgebras, specific to this characteristic, constructed in terms
of orthogonal and symplectic triple systems, which are defined in
terms of central simple degree three Jordan algebras.
\end{abstract}

\maketitle


\section{Introduction}

In \cite{Eld3}, some simple finite dimensional Lie superalgebras
over fields of characteristic $3$, with no counterpart in Kac's
classification (\cite{Kac77}), were constructed by means of the so
called symplectic and orthogonal triple systems over these fields,
most of them related to simple Jordan algebras.

On the other hand, Freudenthal Magic Square, which contains in
characteristic $0$ the exceptional simple finite dimensional Lie
algebras, other than $G_2$, is usually constructed based on two
ingredients: a unital composition algebra and a central simple
degree $3$ Jordan algebra (see \cite[Chapter IV]{Schafer}). This
construction, due to Tits, does not work in characteristic $3$.

A more symmetric construction, based on two unital composition
algebras, which play symmetric roles, and their triality Lie
algebras, has been given recently by several authors (\cite{AF93},
\cite{BS1}, \cite{LM1}). Among other things, this construction has
the advantage of being valid too in characteristic $3$. Simpler
formulas for triality appear if symmetric composition algebras are
used, instead of the more classical unital composition algebras.
These algebras permit a simple construction of Freudenthal Magic
Square (\cite{EldIbero1}).

But the characteristic $3$ presents an exceptional feature, as only
over fields of this characteristic there are nontrivial composition
superalgebras, which appear in dimensions $3$ and $6$. This fact
allows to extend Freudenthal Magic Square (\cite{CunEld1}) with the
addition of two further rows and columns, filled with (mostly
simple) Lie superalgebras.

This paper is devoted to show that the Lie superalgebras (or their
derived subalgebras) that appear in this extended Freudenthal Magic
Square and which are constructed in terms of a three dimensional
composition superalgebra and a composition algebra are among the
simple Lie superalgebras defined in \cite{Eld3} by means of simple
orthogonal triple systems; while those that appear constructed in
terms of a six dimensional composition superalgebra and a
composition algebra are among those in \cite{Eld3} defined by means
of simple symplectic triple systems. Some reflections will be given
too on the simple Lie superalgebras obtained from two nontrivial
composition superalgebras.

\smallskip

The paper is structured as follows. Section 2 will be devoted to
review the construction of the extended Freudenthal Magic Square in
characteristic $3$ in \cite{CunEld1}. Then Section 3 will deal with
the Lie superalgebras of derivations of the Jordan superalgebras of
$3\times 3$ hermitian matrices over unital composition
superalgebras, and their relationship to the first row of the
extended Freudenthal Magic Square. Section 4 will will be devoted to
orthogonal triple systems, some of the simple Lie superalgebras
constructed out of them in \cite{Eld3}, and their connections to the
column of the extended Freudenthal Magic Square which correspond to
the three dimensional composition superalgebras. Section 5 deals
with symplectic triple systems and the column attached to the six
dimensional composition superalgebras. Finally, in Section 6, the
relationship of the remaining Lie superalgebras in the extended
Freudenthal Magic Square to Lie superalgebras in \cite{Eld3} and
\cite{EldPac} will be highlighted, as well as their connections to
triple systems of a mixed nature: the orthosymplectic triple
systems, which will be defined here.


\section{The extended Freudenthal Magic Square in characteristic
$3$}

A quadratic superform on a $\bZ_2$-graded vector space
$U=U\subo\oplus U\subuno$ over a field $k$ is a pair
$q=(q\subo,\bil)$ where $q\subo :U\subo\rightarrow k$ is a quadratic
form, and
 $\bil:U\times U\rightarrow k$ is a supersymmetric even bilinear form
such that $\bil\vert_{U\subo\times U\subo}$ is the polar of $q\subo$
($\bil(x\subo,y\subo)=q\subo(x\subo+y\subo)-q\subo(x\subo)-q\subo(y\subo)$
for any $x\subo,y\subo\in U\subo$).

The quadratic superform $q=(q\subo,\bil)$ is said to be
\emph{regular} if $q\subo$ is regular (definition as in
\cite[p.~xix]{KMRT}) and the restriction of $\bil$ to $U\subuno$ is
nondegenerate.

\smallskip

Then a superalgebra $C=C\subo\oplus C\subuno$ over $k$, endowed with
a regular quadratic superform $q=(q\subo,\bil)$, called the
\emph{norm}, is said to be a \emph{composition superalgebra} (see
\cite{EldOkuCompoSuper}) in case
\begin{subequations}\label{eq:norm}
\begin{align}
&q\subo(x\subo y\subo)=q\subo(x\subo)q\subo(y\subo),\label{eq:qcompo1}\\
&\bil(x\subo y,x\subo z)=q\subo(x\subo)\bil(y,z)=\bil(yx\subo,zx\subo),\label{eq:qcompo2}\\
&\bil(xy,zt)+(-1)^{\lvert x\rvert\lvert y\rvert +\lvert
x\rvert\lvert
 z\rvert+\lvert y\rvert\lvert z\rvert}\bil(zy,xt)=(-1)^{\lvert
 y\rvert\lvert z\rvert}\bil(x,z)\bil(y,t),\label{eq:qcompo3}
\end{align}
\end{subequations}
for any $x\subo,y\subo\in C\subo$ and homogeneous elements
$x,y,z,t\in C$ (where $\lvert x\rvert$ denotes the parity of the
homogeneous element $x$).

The unital composition superalgebras are termed \emph{Hurwitz
superalgebras}, and a composition superalgebra is said to be
\emph{symmetric} in case its bilinear form is associative, that is,
$ \bil(xy,z)=\bil(x,yz), $ for any $x,y,z$.

\smallskip

Hurwitz algebras are the well-known algebras that generalize the
classical real division algebras of the real and complex numbers,
quaternions and octonions. Over any algebraically closed field $k$,
there are exactly four of them: $k$, $k\times k$, $\Mat_2(k)$ and
$C(k)$ (the split Cayley algebra), with dimensions $1$, $2$, $4$ and
$8$.

\smallskip

Only over fields of characteristic $3$  there appear nontrivial
Hurwitz superalgebras (see \cite{EldOkuCompoSuper}):

\begin{itemize}

\item Let $V$ be a two dimensional vector space over a field $k$,
endowed with a nonzero alternating bilinear form $\langle .\vert
.\rangle$.  Consider the superspace $B(1,2)$ (see \cite{She97}) with
\begin{equation}\label{eq:B12a}
B(1,2)\subo =k1,\qquad\text{and}\qquad B(1,2)\subuno= V,
\end{equation}
endowed with the supercommutative multiplication given by
\[
1x=x1=x\qquad\text{and}\qquad uv=\langle u\vert v\rangle 1
\]
for any $x\in B(1,2)$ and $u,v\in V$, and with the quadratic
superform $q=(q\subo,\bil)$ given by:
\begin{equation}\label{eq:B12b}
q\subo(1)=1,\quad \bil(u,v)=\langle u\vert v\rangle,
\end{equation}
for any $u,v\in V$. If the characteristic of $k$ is $3$, then
$B(1,2)$ is a Hurwitz superalgebra (\cite[Proposition
2.7]{EldOkuCompoSuper}).

\smallskip

\item Moreover, with $V$ as before, let $f\mapsto \bar f$ be the
associated symplectic involution on $\End_k(V)$ (so $\langle
f(u)\vert v\rangle =\langle u\vert\bar f(v)\rangle$ for any $u,v\in
V$ and $f\in\End_k(V)$). Consider the superspace $B(4,2)$ (see
\cite{She97}) with
\begin{equation}\label{eq:B42}
B(4,2)\subo=\End_k(V),\qquad\text{and}\qquad B(4,2)\subuno=V,
\end{equation}
with multiplication given by the usual one (composition of maps) in
$\End_k(V)$, and by
\[
\begin{split}
&v\cdot f=f(v)=\bar f\cdot v,\\
&u\cdot v=\langle .\vert u\rangle v\,(w\mapsto \langle w\vert
u\rangle v)\,\in\End_k(V),
\end{split}
\]
for any $f\in\End_k(V)$ and $u,v\in V$; and with quadratic superform
$q=(q\subo,\bil)$ such that
\[
q\subo(f)=\det(f),\qquad\bil(u,v)=\langle u\vert v\rangle,
\]
for any $f\in \End_k(V)$ and $u,v\in V$. Again, if the
characteristic is $3$, $B(4,2)$ is a Hurwitz superalgebra
(\cite[Proposition 2.7]{EldOkuCompoSuper}).

\end{itemize}

\smallskip

Given any Hurwitz superalgebra $C$ with norm $q=(q\subo,\bil)$, its
standard involution is given by $x\mapsto \bar x=\bil(x,1)1-x$. If
$\varphi$ is any automorphism of $C$ with $\varphi^3=1$, then a new
product can be defined on $C$ by means of
\begin{equation}\label{eq:Petersson}
x\bullet y=\varphi(\bar x)\varphi^2(\bar y).
\end{equation}
The resulting superalgebra, denoted by $\bar C_\varphi$, is called a
\emph{Petersson superalgebra} and it turns out to be a symmetric
composition superalgebra.

In particular, for $\varphi=1$, $\bar C=\bar C_1$ is said to be the
\emph{para-Hurwitz superalgebra} attached to $C$.

Over a field $k$ of characteristic $3$, consider the Hurwitz
superalgebra $B(1,2)$, and take a basis $\{v,w\}$ of $V$ with
$\langle v\vert w\rangle =1$. Then, for any scalar $\lambda$, the
even linear map $\varphi:B(1,2)\rightarrow B(1,2)$ such that
$\varphi(1)=1$, $\varphi(v)=v$ and $\varphi(w)=\lambda v+w$, is an
order $3$ automorphism (or $\varphi=1$ if $\lambda =0$). Denote by
$S_{1,2}^\lambda$ the symmetric composition superalgebra
$\overline{B(1,2)}_\varphi$. Also, denote by $S_{4,2}$ the
para-Hurwitz superalgebra $\overline{B(4,2)}$.

Any symmetric composition algebra with a nonzero idempotent is a
Petersson superalgebra (see \cite{EldChema96}, \cite[Chapter
VIII]{KMRT} and \cite{EldOkuCompoSuper}) and this is always the case
over algebraically closed fields. If $S_r$ ($r=1$, $2$, $4$ or $8$)
denotes the para-Hurwitz algebra attached to the split Hurwitz
algebra of dimension $r$ (which is either $k$, $k\times k$,
$\Mat_2(k)$ or $C(k)$), and $\tilde S_8$ denotes the split Okubo
algebra (that is, the pseudo-octonion algebra $P_8(k)$ in
\cite{EldChema96}), which is a particular instance of Petersson
algebra constructed from $C(k)$, then (see \cite[Theorem
4.3]{EldOkuCompoSuper}):

\begin{theorem}\label{th:S1S42} Let $k$ be an algebraically closed
field of characteristic $3$. Then, up to isomorphism, any symmetric
composition superalgebra is one of $S_1$, $S_2$, $S_4$, $S_8$,
$\tilde S_8$, $S_{1,2}^\lambda$ ($\lambda\in k$), or $S_{4,2}$.
\end{theorem}

\smallskip

Given a symmetric composition superalgebra $S$, its \emph{triality
Lie superalgebra} $\tri(S)=\tri(S)\subo\oplus\tri(S)\subuno$ is
defined by:
\begin{multline*}
\tri(S)_{\bar i}=\{ (d_0,d_1,d_2)\in\frosp(S,q)^3_{\bar i}:\\
d_0(x\bullet y)=d_1(x)\bullet y+(-1)^{i\lvert x\rvert}x\bullet
d_2(y)\ \forall x,y\in S\subo\cup S\subuno\},
\end{multline*}
where $\bar i= \bar 0,\bar 1$, and $\frosp(S,q)$ denotes the
associated orthosymplectic Lie superalgebra. The bracket in
$\tri(S)$ is given componentwise.

Now, given two symmetric composition superalgebras $S$ and $S'$ over
a field $k$ of characteristic $\ne 2$, one can form (see \cite[\S
3]{CunEld1}) the Lie superalgebra:
\[
\frg=\frg(S,S')=\bigl(\tri(S)\oplus\tri(S')\bigr)\oplus\bigl(\oplus_{i=0}^2
\iota_i(S\otimes S')\bigr),
\]
where $\iota_i(S\otimes S')$ is just a copy of $S\otimes S'$
($i=0,1,2$),  with bracket given by:

\begin{itemize}
\item $\tri(S)\oplus\tri(S')$ is a Lie subsuperalgebra of $\frg$,
\smallskip

\item $[(d_0,d_1,d_2),\iota_i(x\otimes
 x')]=\iota_i\bigl(d_i(x)\otimes x'\bigr)$,
\smallskip

\item
 $[(d_0',d_1',d_2'),\iota_i(x\otimes
 x')]=(-1)^{\lvert d_i'\rvert\lvert x\rvert}\iota_i\bigl(x\otimes d_i'(x')\bigr)$,
\smallskip

\item $[\iota_i(x\otimes x'),\iota_{i+1}(y\otimes y')]=(-1)^{\lvert
x'\rvert\lvert y\rvert}
 \iota_{i+2}\bigl((x\bullet y)\otimes (x'\bullet y')\bigr)$ (indices modulo
 $3$),
\smallskip

\item $[\iota_i(x\otimes x'),\iota_i(y\otimes y')]=
 (-1)^{\lvert x\rvert\lvert x'\rvert+\lvert x\rvert\lvert y'\rvert +
 \lvert y\rvert\lvert y'\rvert}
 \bil'(x',y')\theta^i(t_{x,y})$ \newline \null\qquad\qquad\qquad\qquad $+
 (-1)^{\lvert y\rvert\lvert x'\rvert}
 \bil(x,y)\theta'^i(t'_{x',y'})$,

\end{itemize}
for any $i=0,1,2$ and homogeneous $x,y\in S$, $x',y'\in S'$,
$(d_0,d_1,d_2)\in\tri(S)$, and $(d_0',d_1',d_2')\in\tri(S')$. Here
$\theta$ denotes the natural automorphism
$\theta:(d_0,d_1,d_2)\mapsto (d_2,d_0,d_1)$ in $\tri(S)$, $\theta'$
the analogous automorphism of $\tri(S')$, and
\begin{equation}\label{eq:txy}
t_{x,y}=\bigl(\sigma_{x,y},\tfrac{1}{2}\bil(x,y)1-r_xl_y,\tfrac{1}{2}\bil(x,y)1-l_xr_y\bigr)
\end{equation}
(with $l_x(y)=x\bullet y$, $r_x(y)=(-1)^{\lvert x\rvert\lvert y
\rvert}y\bullet x$, $\sigma_{x,y}(z)=(-1)^{\lvert y\rvert\lvert
z\rvert}\bil(x,z)y-(-1)^{\lvert x\rvert(\lvert y\rvert +\lvert
z\rvert)}\bil(y,z)x$ for homogeneous $x,y,z\in S$), while $\theta'$
and $t'_{x',y'}$ denote the analogous elements for $\tri(S')$.

\smallskip

This construction is a superization of the algebra construction in
\cite{EldIbero1}, which in turn is based on previous constructions
(with Hurwitz algebras) in \cite{BS1,BS2,LM1,LM2}. It gives a
symmetric and simple construction of Freudenthal Magic Square. The
advantage of using symmetric composition (super)algebras lies in the
simplicity of the formulas needed.

Over an algebraically closed field $k$ of characteristic $3$, and
because of \cite[Theorem 12.2]{EldANewLook}, it is enough to deal
with the Lie superalgebras $\frg(S,S')$, where $S$ and $S'$ are one
of $S_1$, $S_2$, $S_4$, $S_8$, $S_{1,2}=S_{1,2}^0$ or $S_{4,2}$ in
Theorem \ref{th:S1S42}. These are displayed in Table
\ref{ta:supersquare}, which has been obtained in \cite{CunEld1}.

\begin{table}[h!]
$$
\vbox{\offinterlineskip
 \halign{\hfil\ $#$\ \hfil&%
 \vreglon #%
 &\hfil\ $#$\ \hfil&\hfil\ $#$\ \hfil
 &\hfil\ $#$\ \hfil&\hfil\ $#$\ \hfil&%
 \vrule  depth 4pt width .5pt #%
 &\hfil\ $#$\ \hfil&\hfil\ $#$\ \hfil\cr
 \bigstrut &width 0pt&S_1&S_2&S_4&S_8&\omit%
    \vrule height 8pt depth 4pt width .5pt&S_{1,2}&S_{4,2}\cr
 &\multispan8{\hreglonfill}\cr
 S_1&&\frsl_2&\frpgl_3&\frsp_6&\frf_4&&\frpsl_{2,2}&\frsp_6\oplus (14)\cr
 \bigstrut S_2&& &\omit$\frpgl_3\oplus \frpgl_3$&\frpgl_6&\tilde \fre_6&%
     &\bigl(\frpgl_3\oplus\frsl_2\bigr)\oplus\bigl(\frpsl_3\otimes (2)\bigr)&
    \frpgl_6\oplus (20)\cr
 \bigstrut S_4&& & &\frso_{12}&\fre_7&
   &\bigl(\frsp_6\oplus\frsl_2\bigr)\oplus\bigl((13)\otimes (2)\bigr)
    &\frso_{12}\oplus spin_{12}\cr
 \bigstrut S_8&& & & &\fre_8&
    &\bigl(\frf_4\oplus\frsl_2\bigr)\oplus\bigl((25)\otimes (2)\bigr)&
      \fre_7\oplus (56)\cr
 \multispan9{\hregletafill}\cr
 \bigstrut S_{1,2}&& & & & & & \frso_7\oplus 2spin_7 &\frsp_8\oplus(40)\cr
 \bigstrut S_{4,2}&& & & & & & & \frso_{13}\oplus spin_{13}\cr}}
$$
\bigskip
\caption{Freudenthal Magic Supersquare (characteristic
$3$)}\label{ta:supersquare}
\end{table}

Since the construction of $\frg(S,S')$ is symmetric, only the
entries above the diagonal are needed. In Table
\ref{ta:supersquare}, $\frf_4,\fre_6,\fre_7,\fre_8$ denote the
simple exceptional classical Lie algebras, $\tilde\fre_6$ denotes a
$78$ dimensional Lie algebras whose derived Lie algebra is the ($77$
dimensional) simple Lie algebra $\fre_6$ (the characteristic is
$3$!). The even and odd parts of the nontrivial superalgebras in the
table which have no counterpart in Kac's classification in
caracteristic $0$ (\cite{Kac77}) are displayed, $spin$ denotes the
spin module for the corresponding orthogonal Lie algebra, while
$(n)$ denotes a module of dimension $n$. Thus, for example,
$\frg(S_4,S_{1,2})$ is a Lie superalgebra whose even part is
(isomorphic to) the direct sum of the symplectic Lie algebra
$\frsp_6(k)$ and of $\frsl_2(k)$, while its odd part is the tensor
product of a $13$ dimensional module for $\frsp_6(k)$ and the
natural $2$ dimensional module for $\frsl_2(k)$.

A precise description of these modules and of the Lie superalgebras
as contragredient Lie superalgebras is given in \cite{CunEld1}.

\smallskip

The main purpose of this paper is to show the relationships of these
superalgebras $\frg(S,S_{1,2})$ and $\frg(S,S_{4,2})$ to some
superalgebras constructed in \cite{Eld3} by means of orthogonal and
symplectic triple systems, and strongly related to some simple
Jordan algebras.


\section{Some Jordan superalgebras and their derivations}

Given any Hurwitz superalgebra $C$ over a field $k$ of
characteristic $\ne 2$, with norm $q=(q\subo,\bil)$ and standard
involution $x\mapsto \bar x$, the superalgebra $H_3(C,*)$ of
$3\times 3$ hermitian matrices over $C$, where $(a_{ij})^*=(\bar
a_{ji})$, is a Jordan superalgebra under the symmetrized product
\begin{equation}\label{eq:Jproduct}
x\circ y= \frac{1}{2}\bigl( xy+(-1)^{\lvert x\rvert\lvert
y\rvert}yx\bigr).
\end{equation}

Let us consider the associated para-Hurwitz superalgebra $S=\bar C$,
with multiplication $a\bullet b=\bar a\bar b$ for any $a,b\in C$.
Then,
\[
\begin{split}
J=H_3(C,*)&=\left\{ \begin{pmatrix} \alpha_0 &\bar a_2& a_1\\
  a_2&\alpha_1&\bar a_0\\ \bar a_1&a_0&\alpha_2\end{pmatrix} :
  \alpha_0,\alpha_1,\alpha_2\in k,\ a_0,a_1,a_2\in S\right\}\\[6pt]
 &=\bigl(\oplus_{i=0}^2 ke_i\bigr)\oplus
     \bigl(\oplus_{i=0}^2\iota_i(S)\bigr),
\end{split}
\]
where
\begin{equation}\label{eq:eisiotas}
\begin{aligned}
e_0&= \begin{pmatrix} 1&0&0\\ 0&0&0\\ 0&0&0\end{pmatrix}, &
 e_1&=\begin{pmatrix} 0&0&0\\ 0&1&0\\ 0&0&0\end{pmatrix}, &
 e_2&= \begin{pmatrix} 0&0&0\\ 0&0&0\\ 0&0&1\end{pmatrix}, \\
 \iota_0(a)&=2\begin{pmatrix} 0&0&0\\ 0&0&\bar a\\
 0&a&0\end{pmatrix},&
 \iota_1(a)&=2\begin{pmatrix} 0&0&a\\ 0&0&0\\
 \bar a&0&0\end{pmatrix},&
 \iota_2(a)&=2\begin{pmatrix} 0&\bar a&0\\ a&0&0\\
 0&0&0\end{pmatrix},
\end{aligned}
\end{equation}
for any $a\in S$. Identify $ke_0\oplus ke_1\oplus ke_2$ to $k^3$ by
means of $\alpha_0e_0+\alpha_1e_1+\alpha_2e_2\simeq
(\alpha_0,\alpha_1,\alpha_2)$. Then the supercommutative
multiplication \eqref{eq:Jproduct} becomes:
\begin{equation}\label{eq:Jniceproduct}
\left\{\begin{aligned}
 &(\alpha_0,\alpha_1,\alpha_2)\circ(\beta_1,\beta_2,\beta_3)=
    (\alpha_0\beta_0,\alpha_1\beta_1,\alpha_2\beta_2),\\
 &(\alpha_0,\alpha_1,\alpha_2)\circ \iota_i(a)
  =\frac{1}{2}(\alpha_{i+1}+\alpha_{i+2})\iota_i(a),\\
 &\iota_i(a)\circ\iota_{i+1}(b)=\iota_{i+2}(a\bullet b),\\
 &\iota_i(a)\circ\iota_i(b)=2\bil(a,b)\bigl(e_{i+1}+e_{i+2}\bigr),
\end{aligned}\right.
\end{equation}
for any $\alpha_i,\beta_i\in k$, $a,b\in S$, $i=0,1,2$, and where
indices are taken modulo $3$.

\smallskip

The aim of this section is to show that the Lie superalgebra of
derivations of $J$ is naturally isomorphic to the Lie superalgebra
$\frg(S_1,S)$.

This is well-known for algebras, as $\frg(S_1,S)$ is isomorphic to
the Lie algebra $\calT(k,H_3(C,*))$ obtained by means of Tits
construction (see \cite{EldIbero1} and \cite{BS2}), and this latter
algebra is, by its own construction, the derivation algebra of
$H_3(C,*)$. What will be done in this section is to make explicit
this isomorphism $\frg(S_1,S)\cong \der J$ and extend it to
superalgebras.

To begin with, \eqref{eq:Jniceproduct} shows that $J$ is graded over
$\bZ_2\times \bZ_2$ with:
\[
J_{(0,0)}=k^3,\quad J_{(1,0)}=\iota_0(S),\quad
 J_{(0,1)}=\iota_1(S),\quad J_{(1,1)}=\iota_2(S)
\]
and, therefore, $\der J$ is accordingly graded over
$\bZ_2\times\bZ_2$:
\[
(\der J)_{(i,j)}=\{ d\in\der J: d\bigl(J_{(r,s)}\bigr)\subseteq
J_{(i+r,j+s)}\ \forall r,s=0,1\}.
\]

\begin{lemma}\label{le:derJ00} $(\der J)_{(0,0)}=\{ d\in \der J : d(e_i)=0\ \forall
i=0,1,2\}$.
\end{lemma}
\begin{proof}
If $d\in(\der J)_{(0,0)}$, then $d(e_i)\in J_{(0,0)}$ for any $i$.
But $J_{(0,0)}$ is isomorphic to $k^3$, whose Lie algebra of
derivations is trivial. Hence, $d(e_i)=0$ for any $i$.

Conversely, if $d\in \der J$ and $d(e_i)=0$ for any $i$, then for
any $i=0,1,2$ and any $s\in S=\bar C$:
\[
\begin{split}
&d\bigl(\iota_i(s)\bigr)=2d\bigl(e_{i+1}\circ \iota_i(s)\bigr)=
  2e_{i+1}\circ d\bigl(\iota_i(s)\bigr),\\
&d\bigl(\iota_i(s)\bigr)=2d\bigl(e_{i+2}\circ \iota_i(s)\bigr)=
  2e_{i+2}\circ d\bigl(\iota_i(s)\bigr),\\
  &0=d\bigl(e_i\circ \iota_i(s)\bigr)=e_i\circ
  d\bigl(\iota_i(s)\bigr),
\end{split}
\]
so $d\bigl(\iota_i(s)\bigr)\in \{x\in J : e_i\circ x=0,\
e_{i+1}\circ x=\frac{1}{2}x=e_{i+2}\circ x\}=\iota_i(S)$. Hence $d$
preserves the grading: $d\in(\der J)_{(0,0)}$.
\end{proof}

\smallskip

\begin{lemma}\label{le:derJ00triS} $(\der J)_{(0,0)}$ is isomorphic to the triality Lie
superalgebra $\tri(S)$.
\end{lemma}
\begin{proof} For any homogeneous $d\in(\der J)_{(0,0)}$, there are
homogeneous linear maps $d_i\in \End_k(S)$, $i=0,1,2$, of the same
parity, such that
$d\bigl(\iota_i(s)\bigr)=\iota_i\bigl(d_i(s)\bigr)$ for any $s\in S$
and $i=0,1,2$. Now, for any $a,b\in S$ and $i=0,1,2$:
\[
\begin{split}
0&=2\bil(a,b)d(e_{i+1}+e_{i+2})=d\bigl(\iota_i(a)\circ\iota_i(b)\bigr)\\
 &=\iota_i\bigl(d_i(a)\bigr)\circ\iota_i(b)+(-1)^{\lvert
 a\rvert\lvert d\rvert}\iota_i(a)\circ\iota_i\bigl(d_i(b)\bigr)\\
 &=2\bigl(\bil(d_i(a),b)+(-1)^{\lvert
   a\rvert\lvert d\rvert}\bil(a,d_i(b)\bigr)(e_{i+1}+e_{i+2}),
\end{split}
\]
so $d_i$ belongs to the orthosymplectic Lie superalgebra
$\frosp(S,\bil)$. Also,
\[
\begin{split}
\iota_i\bigl(d_i(a\bullet b)\bigr)&=
    d\bigl(\iota_i(a\bullet b)\bigr)=d\bigl(\iota_{i+1}(a)\circ
    \iota_{i+2}(b)\bigr)\\
  &=d\bigl(\iota_{i+1}(a)\bigr)\circ \iota_{i+2}(b) +
    (-1)^{\lvert a\rvert\lvert d\rvert}\iota_{i+1}(a)\circ
    d\bigl(\iota_{i+2}(b)\bigr)\\
  &=\iota_{i+1}\bigl(d_{i+1}(a)\bigr)\circ \iota_{i+2}(b) +
    (-1)^{\lvert a\rvert\lvert d\rvert}\iota_{i+1}(a)\circ
    \iota_{i+2}\bigl(d_{i+2}(b)\bigr)\\
  &=\iota_i\Bigl(d_{i+1}(a)\bullet b+(-1)^{\lvert
       a\rvert\lvert d\rvert}a\bullet d_{i+2}(b)\Bigr),
\end{split}
\]
which shows that $(d_0,d_1,d_2)\in \tri(S)$. Now, the linear map
\[
\begin{split}
\tri(S)&\longrightarrow (\der J)_{(0,0)}\\
 (d_0,d_1,d_2)&\mapsto D_{(d_0,d_1,d_2)},
\end{split}
\]
such that
\begin{equation}\label{eq:Dd0d1d2}
\left\{\begin{aligned} &D_{(d_0,d_1,d_2)}(e_i)=0,\\
   &D_{(d_0,d_1,d_2)}\bigl(\iota_i(a)\bigr)=
   \iota_i\bigl(d_i(a)\bigr)
   \end{aligned}\right.
\end{equation}
for any $i=0,1,2$ and $a\in S$, is clearly an isomorphism.
\end{proof}

\smallskip

For any $i=0,1,2$ and $a\in S$, consider the following inner
derivation of the Jordan superalgebra $J$:
\begin{equation}\label{eq:Dia}
D_i(a)=2\bigl[ L_{\iota_i(a)},L_{e_{i+1}}\bigr]
\end{equation}
(indices modulo $3$), where $L_x$ denotes the multiplication by $x$
in $J$. Note that the restriction of $L_{e_i}$ to
$\iota_{i+1}(S)\oplus\iota_{i+2}(S)$ is half the identity, so the
inner derivation $\bigl[ L_{\iota_i(a)},L_{e_{i}}\bigr]$ is trivial
on $\iota_{i+1}(S)\oplus\iota_{i+2}(S)$, which generates $J$. Hence
\begin{equation}\label{eq:Liotaiaei}
\bigl[ L_{\iota_i(a)},L_{e_{i}}\bigr]=0
\end{equation}
for any $i=0,1,2$ and $a\in S$. Also, $L_{e_0+e_1+e_2}$ is the
identity map, so $\bigl[ L_{\iota_i(a)},L_{e_0+e_{1}+e_2}\bigr]=0$,
and hence
\begin{equation}\label{eq:Liotaiaei2}
D_i(a)=2\bigl[ L_{\iota_i(a)},L_{e_{i+1}}\bigr]=-2\bigl[
L_{\iota_i(a)},L_{e_{i+2}}\bigr].
\end{equation}

A straightforward computation with \eqref{eq:Jniceproduct} gives
\begin{equation}\label{eq:Diaaction}
\begin{split}
&D_i(a)(e_i)=0,\ D_i(a)(e_{i+1})=\frac{1}{2} \iota_i(a),\
  D_i(a)(e_{i+2})=-\frac{1}{2}\iota_i(a),\\
&D_i(a)\bigl(\iota_{i+1}(b)\bigr)=-\iota_{i+2}(a\bullet b),\\
&D_i(a)\bigl(\iota_{i+2}(b)\bigr)=(-1)^{\lvert a\rvert\lvert
  b\rvert}\iota_{i+1}(b\bullet a),\\
&D_i(a)\bigl(\iota_i(b)\bigr)=2\bil(a,b)(-e_{i+1}+e_{i+2}),
\end{split}
\end{equation}
for any $i=0,1,2$ and any homogeneous elements $a,b\in S$.

Denote by $D_i(S)$ the linear span of the $D_i(a)$'s, $a\in S$.

\begin{lemma}\label{le:D0D1D2S}
$D_0(S)=(\der J)_{(1,0)}$, $D_1(S)=(\der J)_{(0,1)}$, and
$D_2(S)=(\der J)_{(1,1)}$.
\end{lemma}
\begin{proof}
By symmetry, it is enough to prove the first assertion. The
containment $D_0(S)\subseteq (\der J)_{(1,0)}$ is clear.

Note that $e_0+e_1+e_2$ is the unity element of $J$, and $e_i\circ
e_i=e_i$ for any $i$. Thus any $d\in \der J$ satisfies
\[
d(e_0+e_1+e_2)=0,\quad\text{and}\quad d(e_i)=2e_i\circ d(e_i),
\]
for any $i=0,1,2$. Hence,
\[
d(e_i)\in \{ x\in J : e_i\circ
x=\frac{1}{2}x\}=\iota_{i+1}(S)\oplus\iota_{i+2}(S).
\]
For $d\in (\der J)_{(1,0)}$ one gets
\[
d(e_0)\in J_{(1,0)}\cap\bigl(\iota_1(S)\oplus \iota_2(S)\bigr)=
 \iota_0(S)\cap\bigl(\iota_1(S)\oplus \iota_2(S)\bigr)=0,
\]
and since $d(e_0+e_1+e_2)=0$, $d(e_1)=-d(e_2)$ follows. Since
$d\in(\der J)_{(1,0)}$, there exists an element $a\in S$ such that
$d(e_1)=\iota_0(a)$, and then
\[
d-D_0(a)\in (\der J)_{(1,0)}\cap\{f\in\der J: f(e_i)=0\ \forall
i=0,1,2\},
\]
and $d=D_0(a)$ by Lemma \ref{le:derJ00}.
\end{proof}

Therefore, the $\bZ_2\times\bZ_2$-grading of $\der J$ becomes
\begin{equation}\label{eq:derJDS}
\der J=D_{\tri(S)}\oplus\bigl(\oplus_{i=0}^2 D_i(S)\bigr)
\end{equation}

On the other hand, $S_1=k1$, with $1\bullet 1=1$ and $\bil(1,1)=2$,
so $\tri(S_1)=0$ and for the para-Hurwitz superalgebra $S$:
\[
\begin{split}
\frg(S_1,S)&=\tri(S)\oplus\bigl(\oplus_{i=0}^2 \iota_i(S_1\otimes
      S)\bigr)\\
   &=\tri(S)\oplus\bigl(\oplus_{i=0}^2\iota_i(1\otimes S)\bigr).
\end{split}
\]

\begin{theorem}\label{th:gS1SderJ} Let $S$ be a para-Hurwitz
superalgebra over a field of characteristic $\ne 2$ and let $J$ be
the Jordan superalgebra of $3\times 3$ hermitian matrices over the
associated Hurwitz superalgebra. Then the linear map:
\begin{equation}\label{eq:PhigS1SderJ}
\Phi:\frg(S_1,S)\longrightarrow \der J,
\end{equation}
such that
\[
\begin{split}
\Phi\bigl((d_0,d_1,d_2)\bigr)&=D_{(d_0,d_1,d_2)},\\
\Phi\bigl(\iota_i(1\otimes a)\bigr)&=D_i(a),
\end{split}
\]
for any $i=0,1,2$, $a\in S$ and $(d_0,d_1,d_2)\in\tri(S)$, is an
isomorphism of Lie superalgebras.
\end{theorem}

\begin{proof}
Equation \eqref{eq:derJDS} shows that $\Phi$ is an isomorphism of
vector spaces. By symmetry it is enough to check that:
\begin{romanenumerate}

\item $\bigl[ D_{(d_0,d_1,d_2)},D_{(f_0,f_1,f_2)}\bigr]
  =D_{([d_0,f_0],[d_1,f_1],[d_2,f_2])}$,

\item $\bigl[D_{(d_0,d_1,d_2)},D_i(a)\bigr]=D_i\bigl(d_i(a)\bigr)$,

\item $\bigl[D_0(a),D_1(b)]=D_2(a\bullet b)$, and

\item $\bigl[D_0(a),D_0(b)\bigr]=2D_{t_{a,b}}$,

\end{romanenumerate}
for any $(d_0,d_1,d_2),(f_0,f_1,f_2)\in\tri(S)$ and $a,b\in S$,
where $t_{a,b}$ is defined in \eqref{eq:txy}.

Both (i) and (ii) are clear from the definitions in
\eqref{eq:Dd0d1d2} and \eqref{eq:Dia}. Now, for any homogeneous
elements $a,b\in S$,
\[
\begin{split}
[D_0(a),D_1(b)]&=\bigl[D_0(a),2[L_{\iota_1(b)},L_{e_2}]\bigr]\\
  &=2\bigl[L_{D_0(a)(\iota_1(b))},L_{e_2}\bigr] +
     (-1)^{\lvert a\rvert\lvert
     b\rvert}2\bigl[L_{\iota_1(b)},L_{D_0(a)(e_2)}\bigr]\\
  &=-2\bigl[L_{\iota_2(a\bullet b)},L_{e_2}\bigr]-(-1)^{\lvert
  a\rvert\lvert b\rvert}\bigl[L_{\iota_1(b)},L_{\iota_0(a)}\bigr]\\
  &=\bigl[L_{\iota_0(a)},L_{\iota_1(b)}\bigr]\quad\text{by
  \eqref{eq:Liotaiaei}},
\end{split}
\]
which is an inner derivation in $(\der J)_{(1,1)}$. But
\[
\bigl[L_{\iota_0(a)},L_{\iota_1(b)}\bigr](e_0)
 =\frac{1}{2}\iota_0(a)\circ\iota_1(b)=
 \frac{1}{2}\iota_2(a\bullet b)=D_2(a\bullet b)(e_0).
\]
which shows, since $(\der J)_{(1,1)}=D_2(S)$ and
$D_2(a)(e_0)=\frac{1}{2}\iota_2(a)$ for any $a\in S$, that
$\bigl[L_{\iota_0(a)},L_{\iota_1(b)}\bigr]=D_2(a\bullet b)$. Hence
(iii) follows.

Finally,
\[
\begin{split}
[D_0(a),D_0(b)]&=\bigl[D_0(a),2[L_{\iota_0(b)},L_{e_1}]\bigr]\\
 &=2\bigl[L_{D_0(a)(\iota_0(b))},L_{e_1}\bigr]+
  (-1)^{\lvert a\rvert\lvert
  b\rvert}2\bigl[L_{\iota_0(b)},L_{D_0(a)(e_1)}\bigr]\\
 &=4\bil(a,b)[L_{e_0-e_2},L_{e_1}]+(-1)^{\lvert a\rvert\lvert
  b\rvert}[L_{\iota_0(b)},L_{\iota_0(a)}]\\
  &=-[L_{\iota_0(a)},L_{\iota_0(b)}].
\end{split}
\]
But for any homogeneous $a,b,x\in S$:
\[
\begin{split}
\bigl[L_{\iota_0(a)},L_{\iota_0(b)}\bigr]\bigl(\iota_1(x)\bigr)&=
  \iota_0(a)\circ\bigl(\iota_0(b)\circ\iota_1(x)\bigr) -
  (-1)^{\lvert a\rvert\lvert
  b\rvert}\iota_0(b)\circ\bigl(\iota_0(a)\circ\iota_1(x)\bigr)\\
 &=\iota_0(a)\circ\iota_2(b\bullet x)-(-1)^{\lvert a\rvert\lvert
  b\rvert}\iota_0(b)\circ\iota_2(a\bullet x)\\
 &=(-1)^{\lvert a\rvert(\lvert b\rvert+\lvert x\rvert)}
   \iota_1\bigl((b\bullet x)\bullet a\bigr)-
   (-1)^{\lvert b\rvert\lvert x\rvert}\iota_1\bigl((a\bullet
   x)\bullet b\bigr)
\end{split}
\]
(see \eqref{eq:Jniceproduct}). Now, from \eqref{eq:qcompo3} and the
nondegeneracy and associativity of the bilinear form of $S$:
\[
(b\bullet x)\bullet a+(-1)^{\lvert a\rvert\lvert x\rvert +
  \lvert b\rvert\lvert x\rvert +\lvert a\rvert\lvert b\rvert}
  (a\bullet x)\bullet b=(-1)^{\lvert b\rvert\lvert x\rvert}\bil(a,b)x,
\]
so
\[
\bigl[L_{\iota_0(a)},L_{\iota_0(b)}\bigr]\bigl(\iota_1(x)\bigr)=
 \iota_1\bigl(-\bil(a,b)x+2r_al_b(x)\bigr)
\]
and
\[
[D_0(a),D_0(b)]\bigl(\iota_1(x)\bigr)=\iota_1\bigl(\bil(a,b)x-2r_al_b(x)\bigr).
\]

In the same vein, one gets
\[
[D_0(a),D_0(b)]\bigl(\iota_2(x)\bigr)=\iota_2\bigl(\bil(a,b)x-2l_ar_b(x)\bigr),
\]
so the restriction of $[D_0(a),D_0(b)]$ to
$\iota_1(S)\oplus\iota_2(S)$ coincides with the restriction of
$2D_{t_{a,b}}$. Since $\iota_1(S)\oplus\iota_2(S)$ generates the
superalgebra $J$, (iv) follows.
\end{proof}

\smallskip

The Lie superalgebra $[L_J,L_J]$ is the Lie superalgebra $\inder J$
of inner derivations of $J$. The proof above shows that $(\der
J)_{(r,s)}=(\inder J)_{(r,s)}$ for $(r,s)\ne (0,0)$, while
\[
(\inder
J)_{(0,0)}=\sum_{i=0}^2\bigl[L_{\iota_i(S)},L_{\iota_i(S)}\bigr]=
D_{\sum_{i=0}^2\theta^i(t_{S,S})}=\Phi\bigl(\sum_{i=0}^2\theta^i(t_{S,S})\bigr)
\]
(recall that $\theta\bigl((d_0,d_1,d_2)\bigr)=(d_2,d_0,d_1)$ for any
$(d_0,d_1,d_2)\in\tri(S)$).

In characteristic $3$, $\tri(S)=\sum_{i=0}^2\theta^i(t_{S,S})$ if
$\dim S=1,4$ or $8$ (\cite{EldIbero1}), and the same happens with
$\tri(S_{1,2})$ and $\tri(S_{4,2})$, because of \cite[Corollaries
2.12 and 2.23]{CunEld1}, while for $\dim S=2$, $\tri(S)$ has
dimension $2$ and $\sum_{i=0}^2\theta^i(t_{S,S})=t_{S,S}$ has
dimension $1$ (see \cite{EldIbero1}). In characteristic $\ne 3$,
$\tri(S)=\sum_{i=0}^2\theta^i(t_{S,S})$ always holds. Therefore, the
proof above, together with the results in \cite{EldIbero1} and
\cite{CunEld1} gives:

\begin{corollary}
Let $S$ be a para-Hurwitz (super)algebra over a field of
characteristic $\ne 2$, and let $J$  be the Jordan (super)algebra of
$3\times 3$ hermitian matrices over the associated Hurwitz
(super)algebra. Then $\der J$ is a simple Lie (super)algebra that
coincides with $\inder J$ unless the characteristic is $3$ and $\dim
S=2$. In this latter case $\inder J$ coincides with $[\der J,\der
J]$, which is a codimension $1$ simple ideal of $\der J$.
\end{corollary}


\section{Orthogonal triple systems}

Orthogonal triple systems were first introduced in \cite[Section
V]{Oku93}:

\begin{definition}\label{df:OTS}
Let $T$ be a vector space over a field $k$ endowed with a nonzero
symmetric bilinear form $(.\vert.):T\times T\rightarrow k$, and a
triple product $T\times T\times T\rightarrow T$: $(x,y,z)\mapsto
[xyz]$. Then $T$ is said to be an \emph{orthogonal triple system} if
it satisfies the following identitities:
\begin{subequations}\label{eq:OTS}
\begin{align}
&[xxy]=0\label{eq:OTSa}\\
&[xyy]=(x\vert y)y-(y\vert y)x\label{eq:OTSb}\\
&[xy[uvw]]=[[xyu]vw]+[u[xyv]w]+[uv[xyw]]\label{eq:OTSc}\\
&([xyu]\vert v)+(u\vert [xyv])=0\label{eq:OTSd}
\end{align}
\end{subequations}
for any elements $x,y,u,v,w\in T$.
\end{definition}

Equation \eqref{eq:OTSc} shows that $\inder T=\espan{[xy.]: x,y\in
T}$ is a subalgebra (actually an ideal) of the Lie algebra $\der T$
of derivations of $T$, whose elements are called \emph{inner
derivations}. Because of \eqref{eq:OTSb}, if $\dim T\geq 2$, then
$\der T$ is contained in the orthogonal Lie algebra
$\frso\bigl(T,(.\vert .)\bigr)$.

The interesting point about these systems is that they provide a
nice construction of Lie superalgebras (see \cite[Theorem
4.5]{Eld3}):

\begin{theorem}\label{th:OTSsuperLie}
Let $T$ be an orthogonal triple system  and let
$\bigl(V,\langle.\vert.\rangle\bigr)$ be a two dimensional vector
space endowed with a nonzero alternating bilinear form. Let $\frs$
be a Lie subalgebra of $\der J$ containing $\inder T$. Define the
superalgebra $\frg=\frg(T,\frs)=\frgo\oplus\frguno$ with
\[
\begin{cases}
\frgo= \frsp(V)\oplus \frs&\text{(direct sum of ideals),}\\
\frguno=V\otimes T\,,
\end{cases}
\]
and superanticommutative multiplication given by:
\begin{itemize}
\item
$\frgo$ is a Lie subalgebra of $\frg$;
\item
$\frgo$ acts naturally on $\frguno$, that is,
\[
[s,v\otimes x]=s(v)\otimes x,\qquad [d,v\otimes x]=v\otimes d(x),
\]
for any $s\in \frsp(V)$, $d\in\frs$, $v\in V$, and $x\in T$;
\item
for any $u,v\in V$ and $x,y\in T$:
\begin{equation}\label{eq:gToproduct}
[u\otimes x,v\otimes y]=
  -(x\vert y)\gamma_{u,v} +\langle u\vert v\rangle d_{x,y}
\end{equation}
where $\gamma_{u,v}=\langle u\vert .\rangle v+\langle v\vert
.\rangle u$ and $d_{x,y}=[xy.]$.
\end{itemize}
Then $\frg(T,\frs)$ is a Lie superalgebra. Moreover, $\frg(T,\frs)$
is simple if and only if $\frs$ coincides with $\inder T$ and $T$ is
simple.

Conversely, given a  Lie superalgebra $\frg=\frgo\oplus\frguno$ with
\[
\begin{cases}
\frgo=\frsp(V)\oplus \frs &\text{(direct sum of ideals),}\\
\frguno=V\otimes T&\text{(as a module for $\frgo$),}
\end{cases}
\]
where $T$ is a module for $\frs$, by $\frsp(V)$-invariance of the
Lie bracket, equation \eqref{eq:gToproduct} is satisfied for a
symmetric bilinear form $(.\vert .):T\times T\rightarrow k$ and a
skew-symmetric bilinear map $d_{.,.}:T\times T\rightarrow \frs$.
Then, if $(.\vert .)$ is not $0$ and a triple product on $T$ is
defined by means of $[xyz]=d_{x,y}(z)$, $T$ becomes and orthogonal
triple system and the image of $\frs$ in $\frgl(T)$ under the given
representation is a subalgebra of $\der T$ containing $\inder T$.
\end{theorem}

Given an orthogonal triple system $T$, the Lie superalgebra
$\frg(T,\inder T)$ will be denoted simply by $\frg(T)$.

\begin{remark} The bracket \eqref{eq:gToproduct} is not exactly the
one that appears in \cite[Equation (4.6)]{Eld3}, but this latter
multiplied by $-1$. This is just obtained by changing
$\langle.\vert.\rangle$ by its negative in \cite{Eld3}.
\end{remark}

\smallskip

The classification of the simple finite dimensional orthogonal
triple systems over fields of characteristic $0$ appears in
\cite[Theorem 4.7]{Eld3}, based on the classification of the simple
$(-1,-1)$ balanced Freudenthal Kantor triple systems in
\cite[Theorem 4.3]{EKO03}.

In characteristic $3$, there appears at least a new family of simple
orthogonal triple systems (see \cite[Examples 4.20]{Eld3}).
Actually, let $C$ be a Hurwitz algebra over a field $k$ of
characteristic $3$, let $S$ be the associated para-Hurwitz algebra
and let $J=H_3(C,*)$ be the simple Jordan algebra already considered
in Section 3. Let $t$ be the natural trace form on $J$ and let
$t(x,y)=t(x\circ y)$ for any $x,y\in J$. Let $J_0=\{ x\in J:
t(x)=0\}$ be the set of trace zero elements. Since the
characteristic is $3$, $1=e_0+e_1+e_2\in J_0$. Consider the quotient
vector space
\[
T_J^o=J_0/k1,
\]
with triple product given by:
\begin{equation}\label{eq:xyzhat}
[\hat x\hat y\hat z]=\widehat{[L_x,L_y](z)},
\end{equation}
for any $x,y,z\in J_0$, where $\hat x$ denotes the class of the
element $x\in J$ modulo $k1$, and with nondegenerate symmetric
bilinear form given by:
\[
(\hat x\vert \hat y)=t(x\circ y)
\]
for any $x,y\in J_0$. These are well defined maps, and with them
$T_J^o$ becomes a simple orthogonal triple system (\cite[Examples
4.20]{Eld3}), and hence $\frg(T_J^o)$ is a simple Lie superalgebra.

The results in \cite[Theorem 4.7(iv) and Theorem 4.23]{Eld3} show
that:
\begin{itemize}

\item for $\dim C=1$, $\frg(T_J^o)$ is isomorphic to the simple
projective special Lie superalgebra $\frpsl_{2,2}(k)$;

\item for $\dim C=2$, $\frg(T_J^o)$ is a simple Lie superalgebra
of dimension $24$, whose even part is the direct sum of a copy of
$\frsl_2(k)$ and of a form of $\frpsl_3(k)$;

\item for $\dim C=4$, $\frg(T_J^o)$ is a simple Lie superalgebra
of dimension $50$, whose even part is the direct sum of a copy of
$\frsl_2(k)$ and of a form of the symplectic Lie algebra
$\frsp_6(k)$;

\item for $\dim C=8$, $\frg(T_J^o)$ is a simple Lie superalgebra
of dimension $105$, whose even part is the direct sum of a copy of
$\frsl_2(k)$ and of a simple Lie algebra $\frf_4$ of type $F_4$.

\end{itemize}

With the exception of $\frpsl_{2,2}(k)$, none of the above simple
Lie superalgebras have counterparts in Kac's classification in
characteristic $0$ (\cite{Kac77}). For $\dim C=2$ and $k$
algebraically closed, the simple Lie superalgebra $\frg(T_J^o)$ has
recently appeared, in a completely different way, in \cite[4.2
Theorem]{BouLeites}, as the Lie superalgebra denoted there by
$\text{bj}$.

\smallskip

The definition of the triple product \eqref{eq:xyzhat} shows that
the Lie algebra of derivations $\der J$, which leaves invariant the
trace form and annihilates $k1$, embeds naturally in $\frgl(T_J^o)$
and, moreover, that $\inder J=\text{span}\{[L_x,L_y]:$ $x,y\in J\}$
maps, under this embedding, onto $\inder T_J^o$.

Although not needed later on, if the dimension of $C$ is $2$, $4$ or
$8$, it can be shown that $\der T_J^o$ coincides with $\der J$ as
follows. Any $d\in \der T_J^o$ induces an even derivation $\delta$
of $\frg(T_J^o)=\bigl(\frsp(V)\oplus\inder T_J^o\bigr)\oplus
(V\otimes T_J^o)$, such that $\delta(s)=0$ for any $s\in \frsp(V)$,
$\delta(f)=[d,f]$ for any $f\in \inder T_J^o$, and $\delta(v\otimes
x)=v\otimes d(x)$ for any $v\in V$ and $x\in T_J^o$. For $\dim C=4$
(respectively $8$), $\inder T_J^o\simeq \inder J=\der J$ is a form
of $\frsp_6(k)$ (respectively, a simple Lie algebra of type $F_4$),
whose Killing form is nondegenerate (even in this characteristic).
Hence its derivations are all inner, and there exists an inner
derivation $\tilde d\in \inder T_J^o$ such that the restriction
$\delta\vert_{\inder T_J^o}=\ad \tilde d$. Hence $\delta-\ad \tilde
d$ is an even derivation of $\frg(T_J^o)$ which is trivial on the
even part. By irreducibility of the odd part as a module for the
even part, its restriction to the odd part is a scalar, which must
be $0$. Hence $\delta=\ad \tilde d$, which shows that $d=\tilde d$
on $T_J^o$, and hence $\der T_J^o=\inder T_J^o$. If the dimension of
$C$ is $2$, then $\inder T_J^o$ is a form of $\frpsl_3(k)$, whose
derivation algebra is $\frpgl_3(k)$. It follows that $\inder
T_J^o\simeq \inder J$ has codimension at most $1$ in $\der T_J^o$.
Since $\inder J$ is a codimension one ideal in $\der J$ in this
case, one gets that $\der T_J^o$ coincides with (the image under the
natural embedding of) $\der J$.

\smallskip

Under the conditions above, denote simply by $\frg(J)$ the Lie
superalgebra $\frg(T_J^o,\der J)=\bigl(\frsp(V)\oplus \der
J\bigr)\oplus (V\otimes T_J^o)$. If $\dim C=1$, $4$ or $8$, then
$\frg(J)=\frg(T_J^o)$ (recall that this is, by definition, the Lie
superalgebra $\frg(T_J^o,\inder T_J^o)=\frg(T_J^o,\inder J)$), while
for $\dim C=2$, $\frg(T_J^o)=[\frg(J),\frg(J)]$ is a codimension one
ideal in $\frg(J)$.

Now Theorem \ref{th:gS1SderJ} shows that $\der J$ is isomorphic to
the Lie algebra $\frg(S_1,S)$, where $S$ is the para-Hurwitz algebra
attached to $C$. This isomorphism extends to an isomorphism between
the Lie superalgebras $\frg(T_J^o)$ and $\frg(S_{1,2},S)$, as it
will be seen shortly.

The results in \cite[Corollary 2.12]{CunEld1} shows that the
triality Lie superalgebra of $S_{1,2}$ is
\[
\tri(S_{1,2})=\{(d,d,d): d\in \frosp(S_{1,2},\bil)\},
\]
and hence is isomorphic to the orthosymplectic Lie superalgebra
$\frosp(S_{1,2},\bil)$, which is spanned by the operators:
\[
\sigma_{x,y}: z\mapsto (-1)^{\lvert y\rvert\lvert
z\rvert}\bil(x,z)y- (-1)^{\lvert x\rvert(\lvert y\rvert+\lvert
z\rvert)}\bil(y,z)x.
\]
Also, $S_{1,2}=k1\oplus V$, $\bil(1,1)=2\,(=-1)$,
$\bil\vert_{V\times V}=\langle .\vert .\rangle$ is a nonzero
alternating form (see \eqref{eq:B12b}) and $\frsp(V)$ is spanned by
the operators:
\[
\gamma_{u,v}: w\mapsto \langle u\vert w\rangle v+\langle v\vert
w\rangle u.
\]
Then (see \cite[(2.15)-(2.17)]{CunEld1}):
\[
\frosp(S_{1,2},\bil)\subo=\sigma_{V,V}\simeq \frsp(V),\quad
\frosp(S_{1,2},\bil)\subuno=\sigma_{1,V}\simeq V,
\]
where $\sigma_{V,V}$ is identified to $\frsp(V)$ because
$\sigma_{V,V}(1)=0$ and
\begin{equation}\label{eq:sigmasgammas}
\sigma_{u,v}\vert_{V\times V}=-\gamma_{u,v},
\end{equation}
for any $u,v\in V$, while $\sigma_{1,V}$ is identified to $V$ by
means of $\sigma_{1,u}\leftrightarrow u$. Note that
(\cite[(2.16)]{CunEld1}):
\begin{equation}\label{eq:sigma1u}
\sigma_{1,u}(1)=-u,\qquad\sigma_{1,u}(v)=-\langle u\vert v\rangle 1,
\end{equation}
for any $u,v\in V$.

Then, for any para-Hurwitz algebra $S$:
\[
\begin{split}
\frg(S_{1,2},S)&=\bigl(\tri(S_{1,2})\oplus\tri(S)\bigr)
    \oplus\bigl(\oplus_{i=0}^2\iota_i(S_{1,2}\otimes S)\bigr)\\
   &=\bigl(\tri(S_{1,2})\oplus\tri(S)\bigr)
    \oplus\Bigl(\oplus_{i=0}^2\bigl(\iota_i(k1\otimes
    S)\oplus\iota_i(V\otimes S)\bigr)\Bigr)\\
   &=\bigl(\tri(S_{1,2})\oplus\tri(S)\bigr)
    \oplus\bigl(\oplus_{i=0}^2\iota_i(k1\otimes
    S)\bigr)\oplus\bigl(\oplus_{i=0}^2\iota_i(V\otimes
    S)\bigr)\Bigr)\\
   &=\Bigl(\tri(S_{1,2})\subo\oplus\bigl(\tri(S)\oplus
    \bigl(\oplus_{i=0}^2\iota(k1\otimes S)\bigr)\bigr)\Bigr)\\
    &\hspace{2in} \oplus
    \Bigl(\tri(S_{1,2})\subuno\oplus\bigl(\oplus_{i=0}^2\iota_i(V\otimes
    S)\bigr)\Bigr)\\
   &\simeq \bigl(\frsp(V)\oplus\frg(S_1,S)\bigr)\oplus\Bigl(V\oplus
   \bigl(\oplus_{i=0}^2 \iota_i(V\otimes S)\bigr)\Bigr),
\end{split}
\]
where the identifications $\tri(S_{1,2})\subo\simeq
\sigma_{V,V}\simeq \frsp(V)$ and $\tri(S_{1,2})\subuno\simeq
\sigma_{1,V}\simeq V$ are used. Actually,
\[
\begin{split}
\frg(S_{1,2},S)\subo&\simeq
\frsp(V)\oplus\frg(S_1,S)\quad\text{(direct sum of ideals)}\\
\frg(S_{1,2},S)\subuno&\simeq
V\oplus\bigl(\oplus_{i=0}^2\iota_i(V\otimes S)\bigr).
\end{split}
\]

On the other hand,
\[
\frg(J)=(\frsp(V)\oplus\der J)\oplus (V\otimes T_J^o),
\]
and
\[
T_J^o=J_0/k1=
 \hat\iota(k)\oplus\bigl(\oplus_{i=0}^2\hat\iota_i(S)\bigr),
\]
where
\[
\begin{split}
&\hat\iota(1)=\widehat{e_0-e_1}=\widehat{e_1-e_2}=\widehat{e_2-e_0}\quad\text{(as
$\widehat{e_0+e_1+e_2}=0$),}\\
&\hat\iota_i(s)=\widehat{\iota_i(s)}\quad\text{for any $i=0,1,2$ and
$s\in S$}
\end{split}
\]
(see \eqref{eq:eisiotas}).

\begin{theorem}\label{th:gS12SgTJ}
Let $S$ be a para-Hurwitz algebra over a field $k$ of characteristic
$3$ and let $J$ be the Jordan algebra of $3\times 3$ hermitian
matrices over the associated Hurwitz algebra. Then the linear map
\[
\Psi:\frg(S_{1,2},S)\longrightarrow \frg(J)
\]
such that
\begin{itemize}

\item $\Psi\vert_{\frsp(V)}$ is the identity map $\frsp(V)\simeq
\tri(S_{1,2})\subo\rightarrow \frsp(V)$,

\item $\Psi\vert_{\frg(S_1,S)}$ is the isomorphism
$\Phi:\frg(S_1,S)\rightarrow \der J$ in Theorem \ref{th:gS1SderJ},

\item $\Psi\vert_V$ is the `identity' map
$V\simeq\tri(S_{1,2})\subuno\rightarrow V\otimes \hat\iota(k):
u\mapsto u\otimes\hat\iota(1)$,

\item $\Psi\bigl(\iota_i(u\otimes s)\bigr)=u\otimes\hat\iota_i(s)$
for any $u\in V$, $s\in S$ and $i=0,1,2$,

\end{itemize}
is a Lie superalgebra isomorphism.
\end{theorem}

\begin{proof} The proof is obtained by straightforward computations.
First, it is clear that $\Psi$ is a bijective linear map, and that
it maps isomorphically the even part of $\frg(S_{1,2},S)$ onto the
even part of $\frg(J)$, because $\Phi$ is an isomorphism (Theorem
\ref{th:gS1SderJ}). In checking that
$\Psi([g\subo,g\subuno])=[\Psi(g\subo),\Psi(g\subuno)]$ for any
$g_i\in \frg(S_{1,2},S)_i$, $i=\bar 0,\bar 1$, the only nontrivial
instances (up to symmetry) to be checked are the following:

\smallskip

\noindent{$\bullet$}\quad $[\iota_0(1\otimes a),\iota_0(u\otimes
b)]=\bil(a,b)t_{1,u}\simeq \bil(a,b)u$ in $\frg(S_{1,2},S)$ for any
$a,b\in S$ and $u\in V\simeq (S_{1,2})\subuno$; so $
\Psi\bigl([\iota_0(1\otimes a),\iota_0(u\otimes b)]\bigr)
 =\bil(a,b)u\otimes\hat\iota(1)=\bil(a,b)u\otimes\widehat{e_1-e_2},
$
while
\[
\begin{split}
[\Psi(\iota_0(1\otimes a)),\Psi(\iota_0(u\otimes b))]
 &=[\Phi(\iota_0(1\otimes a)),\Phi(\iota_0(u\otimes b))]\\
 &=[D_0(a),u\otimes\hat\iota_0(b)]=u\otimes\widehat{D_0(a)(\iota_0(b))}\\
 &=u\otimes2\bil(a,b)\widehat{(-e_1+e_2)}=\bil(a,b)u\otimes\widehat{e_1-e_2},
\end{split}
\]
because of \eqref{eq:Diaaction}.

\smallskip

\noindent{$\bullet$}\quad $[\iota_0(1\otimes a),\iota_1(u\otimes
b)]=\iota_2(-u\otimes a\bullet b)$ for any $a,b\in S$ and $u\in V$
(as $1\bullet u=\bar 1\bar u=1(-u)=-u$), so
$\Psi\bigl([\iota_0(1\otimes a),\iota_1(u\otimes
b)]\bigr)=-u\otimes\hat\iota_2( a\bullet b)$, while
\[
\begin{split}
[\Psi(\iota_0(1\otimes a)),\Psi(\iota_1(u\otimes b))]
 &=[D_0(a),u\otimes\hat\iota_1(b)]\\
 &=u\otimes\widehat{D_0(a)(\iota_1(b))}\\
 &=-u\otimes \hat\iota_2(a\bullet b)\ \text{(by
 \eqref{eq:Diaaction}).}
\end{split}
\]

\smallskip

\noindent{$\bullet$}\quad Similarly, $\Psi\bigl([\iota_0(1\otimes
a),\iota_2(u\otimes b)]\bigr)=\Psi\bigl(\iota_1(u\otimes b\bullet
a)\bigr)=u\otimes\hat\iota_1( b\bullet a)$, while
$[\Psi(\iota_0(1\otimes a)),\Psi(\iota_2(u\otimes b))]
 =[D_0(a),u\otimes\hat\iota_2(b)]
 =u\otimes \hat\iota_1(b\bullet a)$ too.

\medskip

Also, the instances to be checked for two odd elements are the
following:

\noindent{$\bullet$}\quad $\Psi\bigl([u,\iota_0(v\otimes
a)]\bigr)=-\langle u\vert v\rangle\Psi(\iota_0(1\otimes a))=-\langle
u\vert v\rangle D_0(a)$ for any $u,v\in V$ and $a\in S$, because
$u\simeq \sigma_{1,u}$ acts on $v$ as indicated in
\eqref{eq:sigma1u}. And
\[
\begin{split}
[\Psi(u),\Psi(\iota_0(v\otimes a))]
 &= [u\otimes\hat\iota(1),v\otimes\hat\iota_0(a)]\\
 &=[u\otimes\widehat{e_1-e_2},v\otimes\hat\iota_0(a)]\\
 &=\langle u\vert v\rangle[L_{e_1-e_2},L_{\iota_2(a)}]\\
 &=-2\langle u\vert v\rangle[L_{\iota_0(a)},L_{e_1}]
  =-\langle u\vert v \rangle D_0(a),
\end{split}
\]
because of \eqref{eq:Liotaiaei2}.

\smallskip

\noindent{$\bullet$}\quad For any $u,v\in V$ and $a,b\in S$,
\[
\begin{split}
\Psi\bigl([\iota_0(u\otimes a),\iota_0(v\otimes b)]\bigr)
 &=\Psi\bigl(\bil(a,b)\sigma_{u,v}+\langle u\vert v\rangle
   t_{a,b}\bigr)\\
 &=-\bil(a,b)\gamma_{u,v}+\langle u\vert v\rangle D_{t_{a,b}},
\end{split}
\]
because of \eqref{eq:sigmasgammas}, while
\[
\begin{split}
[\Psi(\iota_0(u\otimes a)),\Psi(\iota_0(v\otimes b))]&=
  [u\otimes\hat\iota_0(a),v\otimes\hat\iota_0(b)]\\
  &=-t(\iota_0(a)\circ\iota_0(b))\gamma_{u,v}+\langle u\vert v\rangle
   [L_{\iota_0(a)},L_{\iota_0(b)}]\\
  &=-\bil(a,b)\gamma_{u,v}-\langle u\vert v\rangle[D_0(a),D_0(b)]\\
  &=-\bil(a,b)\gamma_{u,v}+\langle u\vert v\rangle D_{t_{a,b}},
\end{split}
\]
because of the formulas in the proof of Theorem \ref{th:gS1SderJ}.

\smallskip

\noindent{$\bullet$}\quad Finally, for $u,v\in V$ and $a,b\in S$,
\[
\begin{split}
\Psi\bigl([\iota_0(u\otimes a),\iota_1(v\otimes b)]\bigr)&=
  \Psi\bigl(\iota_2(\langle u\vert v\rangle 1\otimes a\bullet
  b)\bigr)\\
 &=\langle u\vert v\rangle \Phi(\iota_2(1\otimes a\bullet b))\\
 &=\langle u\vert v\rangle D_2(a\bullet b),
\end{split}
\]
while
\[
\begin{split}
[\Psi(\iota_0(u\otimes a)),\Psi(\iota_1(v\otimes b))]
  &=[u\otimes\hat\iota_0(a),v\otimes\hat\iota_1(b)]\\
  &=\langle u\vert v\rangle [L_{\iota_0(a)},L_{\iota_1(b)}]\\
  &= \langle u\vert v\rangle [D_0(a),D_1(b)]\\
  &=\langle u\vert v\rangle D_2(a\bullet b),
\end{split}
\]
again by the arguments in the proof of Theorem \ref{th:gS1SderJ}.
\end{proof}

\medskip

This result immediately shows that the Lie superalgebras
$\frg(S_{1,2},S_r)$ in \cite{CunEld1} are, essentially, the new
simple Lie superalgebras in \cite[Theorem 4.23]{Eld3}. More
specifically:

\begin{corollary} Let $S_r$ ($r=1,2,4,8$) denote the unique
para-Hurwitz superalgebra of dimension $r$ over an algebraically
closed field $k$ of characteristic $3$. Then:
\begin{romanenumerate}
\item $\frg(S_{1,2},S_1)$ is isomorphic to the classical Lie
superalgebra $\frpsl_{2,2}(k)$.

\item $[\frg(S_{1,2},S_2),\frg(S_{1,2},S)]$ is isomorphic to the
simple Lie superalgebra in \cite[Theorem 4.23(i)]{Eld3}, obtained as
$\frg(T_J^o)$ for $J=H_3(k\times k,*)$.

\item $\frg(S_{1,2},S_4)$ is isomorphic to the
simple Lie superalgebra in \cite[Theorem 4.23(iii)]{Eld3}, obtained
as $\frg(T_J^o)$ for $J=H_3(\Mat_2(k),*)$.

\item $\frg(S_{1,2},S_8)$ is isomorphic to the
simple Lie superalgebra in \cite[Theorem 4.23(iv)]{Eld3}, obtained
as $\frg(T_J^o)$ for $J=H_3(C(k),*)$.

\end{romanenumerate}
\end{corollary}


\section{Symplectic triple systems}

Symplectic triple systems appeared for the first time in
\cite{YA75}:

\begin{definition}\label{df:STS}
Let $T$ be a vector space over a field $k$ endowed with a nonzero
alternating bilinear form $(.\vert.):T\times T\rightarrow k$, and a
triple product $T\times T\times T\rightarrow T$: $(x,y,z)\mapsto
[xyz]$. Then $T$ is said to be a \emph{symplectic triple system}
 if it satisfies the following identities:
\begin{subequations}\label{eq:STS}
\begin{align}
&[xyz]=[yxz]\label{eq:STSa}\\
&[xyz]-[xzy]=(x\vert z)y-(x\vert y)z+2(y\vert z)x\label{eq:STSb}\\
&[xy[uvw]]=[[xyu]vw]+[u[xyv]w]+[uv[xyw]]\label{eq:STSc}\\
&([xyu]\vert v)+(u\vert [xyv])=0\label{eq:STSd}
\end{align}
\end{subequations}
for any elements $x,y,z,u,v,w\in T$.
\end{definition}

These systems are strongly related to Freudenthal triple systems and
to Faulkner ternary algebras (see \cite[Theorems 2.16 and
2.18]{Eld3} and references therein).

As for orthogonal triple systems, for any $x,y$ in a symplectic
triple system, the linear map $[xy.]$ is a derivation, and the span
$\inder T$ of these derivations is an ideal of $\der T$, whose
elements are called \emph{inner derivations}. Then (see
\cite[Theorem 2.9]{Eld3}):

\begin{theorem}\label{th:gTSTS}
Let $T$ be a symplectic triple system and let
$\bigl(V,\langle.\vert.\rangle\bigr)$ be a two dimensional vector
space endowed with a nonzero alternating bilinear form. Let $\frs$
be a subalgebra of $\der T$ containing $\inder T$. Define the
$\bZ_2$-graded algebra $\frg=\frg(T,\frs)=\frgo\oplus\frguno$ with
\[
\begin{cases}
\frgo= \frsp(V)\oplus \frs&\text{(direct sum of ideals),}\\
\frguno=V\otimes T\,,
\end{cases}
\]
and anticommutative multiplication given by:
\begin{itemize}
\item
$\frgo$ is a Lie subalgebra of $\frg$,
\item
$\frgo$ acts naturally on $\frguno$; that is
\[
[s,v\otimes x]=s(v)\otimes x,\qquad [d,v\otimes x]=v\otimes d(x),
\]
for any $s\in \frsp(V)$, $d\in\frs$, $v\in V$, and $x\in T$.
\item
For any $u,v\in V$ and $x,y\in T$:
\begin{equation}\label{eq:gTsproduct}
[u\otimes x,v\otimes y]=
  (x\vert y)\gamma_{u,v} +\langle u\vert v\rangle d_{x,y}
\end{equation}
where, as before, $\gamma_{u,v}=\langle u\vert .\rangle v+\langle
v\vert .\rangle u$ and $d_{x,y}=[xy.]$.
\end{itemize}
Then $\frg(T,\frs)$ is a Lie algebra. Moreover, $\frg(T,\frs)$ is
simple if and only if $\frs$ coincides with $\inder T$ and $T$ is
simple.

Conversely, given a $\bZ_2$-graded Lie algebra
$\frg=\frgo\oplus\frguno$ with
\[
\begin{cases}
\frgo=\frsp(V)\oplus \frs &\text{(direct sum of ideals),}\\
\frguno=V\otimes T&\text{(as a module for $\frgo$),}
\end{cases}
\]
where $T$ is a module for $\frs$, by $\frsp(V)$-invariance of the
Lie bracket, equation \eqref{eq:gTsproduct} is satisfied for an
alternating bilinear form $(.\vert .):T\times T\rightarrow k$ and a
symmetric bilinear map $d_{.,.}:T\times T\rightarrow \frs$. Then, if
$(.\vert .)$ is not $0$ and a triple product on $T$ is defined by
means of $[xyz]=d_{x,y}(z)$, $T$ becomes a symplectic triple system,
and the image of $\frs$ in $\frgl(T)$ under the given representation
is a subalgebra of $\der T$ containing $\inder T$.
\end{theorem}

\smallskip

The classification of the simple finite dimensional symplectic
triple systems in characteristic $3$ appears in \cite[Theorem
2.32]{Eld3}, based on the classification of Freudenthal triple
systems in this characteristic by Brown \cite{Bro84}. The symplectic
triple systems which most interest us are the following. Let $C$ be
a Hurwitz algebra over a field $k$ (characteristic $\ne 2$), and let
$J$ be the Jordan algebra of $3\times 3$ hermitian matrices over
$C$, as in Section 2, with the usual trace $t$. Consider the `cross
product' on $J$ define by
\[
x\times y=2x\circ y+t(x)y+t(y)x-s(x,y)1,
\]
where $s(x,y)=t(x)t(y)-t(x\circ y)$. Take the vector space
\begin{equation}\label{eq:TJs}
T_J^s=\left\{\begin{pmatrix} \alpha&a\\ b&\beta\end{pmatrix} :
\alpha,\beta\in k,\, a,b\in J\right\},
\end{equation}
endowed with the alternating bilinear form and triple product such
that, for $x_i=\bigl(\begin{smallmatrix} \alpha_i&a_i\\
b_i&\beta_i\end{smallmatrix}\bigr)\in T_J^s$, $i=1,2,3$:
\[
\begin{aligned}
&(x_1\vert x_2)=
  \alpha_1\beta_2-\alpha_2\beta_1-t(a_1,b_2)+t(b_1,a_2),\\
&[x_1x_2x_3]=\begin{pmatrix} \gamma&c\\ d&\delta
    \end{pmatrix}\quad\text{with}\\
&\quad \gamma=\Bigl(-3(\alpha_1\beta_2+\alpha_2\beta_1)
   +t(a_1,b_2)+t(a_2,b_1)\Bigr)\alpha_3\\
&\qquad\qquad +2\Bigl(\alpha_1t(b_2,a_3)
   +\alpha_2t(b_1,a_3)-t(a_1\times a_2,a_3)\Bigr)\\
&\quad c=\Bigl(-(\alpha_1\beta_2+\alpha_2\beta_1)
   +t(a_1,b_2)+t(a_2,b_1)\Bigr)a_3\\
&\qquad\qquad +2\Bigl(\bigl(t(b_2,a_3)-\beta_2\alpha_3\bigr)a_1+
    \bigl(t(b_1,a_3)-\beta_1\alpha_3\bigr)a_2\Bigr)\\
&\qquad\qquad +2\Bigl( \alpha_1(b_2\times b_3)
   +\alpha_2(b_1\times b_3)+\alpha_3(b_1\times b_2)\Bigr)\\
&\qquad\qquad -2\Bigl( (a_1\times a_2)\times b_3 +
  (a_1\times a_3)\times b_2 +(a_2\times a_3)\times b_1\Bigr)\\
&\quad \text{$\delta=-\gamma^\sigma$, $d=-c^\sigma$, where
$\sigma=(\alpha\beta)(ab)$ (that is, $\gamma^{\sigma}$ and
$c^\sigma$}\\
&\quad\quad\text{are obtained from $\gamma$ and $c$ by
interchanging $\alpha$ and $\beta$ and}\\
&\quad\quad\text{also $a$ and $b$ throughout).} \qed
\end{aligned}
\]

\smallskip

In characteristic $3$ it turns out that any symplectic triple system
is an anti-Lie triple system, and hence the next result
(\cite[Theorem 3.1]{Eld3}) holds:

\begin{theorem}\label{th:STSsuperLie}
Let $T$ be a symplectic triple system over a field $k$ of
characteristic $3$ and let $\frs$ be a subalgebra of $\der T$
containing $\inder T$. Define the superalgebra
$\frtg=\frtg(T,\frs)=\frtg\subo\oplus\frtg\subuno$, with:
\[
\frtg\subo=\frs,\qquad\frtg\subuno=T,
\]
and superanticommutative multiplication given by:
\begin{itemize}
\item
$\frtg\subo$ is a Lie subalgebra of $\frtg$;
\item
$\frtg\subo$ acts naturally on $\frtg\subuno$: $[d,x]=d(x)$ for any
$d\in\frs$ and $x\in T$;
\item
$[x,y]=d_{x,y}=[xy.]$, for any $x,y\in T$.
\end{itemize}
Then $\frtg(T,\frs)$ is a Lie superalgebra. Moreover,
$\frtg(T,\frs)$ is simple if and only if so is $T$ and $\frs=\inder
T$.
\end{theorem}

That is, in Theorem \ref{th:gTSTS} one may ``delete'' $\frsp(V)$ and
$V$ and then the $\bZ_2$-graded Lie algebra there becomes a Lie
superalgebra!

In case $\frs=\inder T$, the Lie superalgebra $\frtg(T,\frs)$ will
just be denoted by $\frtg(T)$.

\smallskip

For the symplectic triple systems $T_J^s$ in \eqref{eq:TJs}, the
results in \cite[Theorem 3.2]{Eld3} show that:

\begin{itemize}

\item If $\dim C=1$, then $\frtg(T_J^s)$ is a simple Lie
superalgebra of dimension $35$, whose even part is a form of the
symplectic Lie algebra $\frsp_6(k)$ and whose odd part is an
irreducible module of dimension $14$ for the even part.

\item If $\dim C=2$, then $\frtg(T_J^s)$ is a simple Lie
superalgebra of dimension $54$, whose even part is a form of the
projective special Lie algebra $\frpsl_6(k)$ and whose odd part is
an irreducible module of dimension $20$ for the even part.

\item If $\dim C=4$, then $\frtg(T_J^s)$ is a simple Lie
superalgebra of dimension $98$, whose even part is a form of the
orthogonal Lie algebra $\frso_{12}(k)$ and whose odd part is an
irreducible module of dimension $32$ (the spin module) for the even
part.

\item If $\dim C=8$, then $\frtg(T_J^s)$ is a simple Lie
superalgebra of dimension $189$, whose even part is a simple
exceptional Lie algebra of type $E_7$ and whose odd part is an
irreducible module of dimension $56$ for the even part.

\end{itemize}

None of the above simple Lie superalgebras have counterparts in
Kac's classification in characteristic $0$.

\smallskip

Since any derivation of $\frsp_6(k)$, $\frso_{12}(k)$ or the simple
Lie algebras of type $E_7$ is inner, it follows, as in Section 4,
that $\inder T_J^s=\der T_J^s$ for $\dim C=1$, $4$ or $8$. For $\dim
C=2$, $\der\frpsl_6(k)=\frpgl_6(k)$. Extending scalars to an
algebraic closure, the symplectic triple system $T_J^s$ is then
isomorphic to the one denoted by $\calT_{\fre_6,\fra_5}$ in \cite[\S
4]{EldIbero2}, which satisfies that $\fra_5=\frpgl_6(k)$ is
contained in its Lie algebra of derivations. It follows that $\der
T_J^s$ is a form of $\frpgl_6(k)$ and $\inder T_J^s=[\der T_J^s,\der
T_J^s]$ is a codimension one ideal in $\der T_J^s$ in this case and,
therefore, $\frtg(T_J^s)=[\frtg(T_J^s,\der T_J^s),\frtg(T_J^s,\der
T_J^s)]$ is a codimension one ideal of $\frtg(T_J^s,\der T_J^s)$.

Assume now that the ground field $k$ is algebraically closed of
characteristic $3$. Then the Lie superalgebras in the extended
Freudenthal Magic Square (Table \ref{ta:supersquare}) have been
constructed in \cite{CunEld1} by means of copies of a two
dimensional vector space $V$, endowed with a nonzero alternating
bilinear form, and copies of the three dimensional simple Lie
algebra $\frsp(V)$, following the approach in \cite{EldIbero2} used
for the Lie algebras in Freudenthal Magic Square.

In particular, if $S_r$ denotes the unique para-Hurwitz algebra of
dimension $r$ over $k$ ($r=1,2,4$ or $8$), the construction of the
Lie algebra $\frg(S_r,S_8)$ uses four copies of $V$ to deal with
$S_8$, say $V_1,V_2,V_3,V_4$. Then (see \cite[\S 4]{EldIbero2}, but
with the indices $1,2,3,4$ for $S_8$ in all cases):
\begin{equation}\label{eq:gSrS8}
\frg(S_r,S_8)=\bigl(\frg(S_r,S_4)\oplus\frsp(V_4)\bigr)\oplus
\Bigl(\bigl(\oplus_{\sigma\in\calS_{S_r,S_8}\setminus
\calS_{S_r,S_4}}\tilde V(\sigma\setminus\{ 4\})\bigr)\otimes
V_4\Bigr),
\end{equation}
where $\calS_{S_8,S_8}=\calS_{\fre_8}$,
$\calS_{S_4,S_8}=\calS_{\fre_7}$, $\calS_{S_2,S_8}=\calS_{\fre_6}$,
and $\calS_{S_1,S_8}=\calS_{\frf_4}$ in \cite{EldIbero2}, and
$\tilde V=V$ unless $r=2$.

According to Theorem \ref{th:gTSTS}, this shows that
$\oplus_{\sigma\in\calS_{S_r,S_8}\setminus \calS_{S_r,S_4}}\tilde
V(\sigma\setminus\{ 4\})$ is a symplectic triple system and
$\frg(S_r,S_4)$ is contained in its Lie algebra of derivations.
These are precisely, up to isomorphism, the simple symplectic triple
systems $T_J^s$ considered so far (see \cite[\S 4]{EldIbero2}).

However, a case by case inspection in \cite[\S 5]{CunEld1} shows
that for any $r=1,2,4,8$,
\[
\calS_{S_r,S_{4,2}}=\calS_{S_r,S_4}\cup\bigl\{\sigma\setminus\{4\}:
\sigma\in\calS_{S_r,S_8}\setminus\calS_{S_r,S_4}\bigr\},
\]
and that, if $\frsp(V_4)$ and $V_4$ are ``deleted'' in
\eqref{eq:gSrS8}, one obtains precisely the Lie superalgebra
$\frg(S_r,S_{4,2})$. These are simple with the exception of $r=2$,
where $[\frg(S_r,S_{4,2}),\frg(S_r,S_{4,2})]$ is a codimension one
simple ideal (see \cite{CunEld1}). Therefore, $\frg(S_r,S_{4,2})$
is, up to isomorphism, the Lie superalgebra $\frtg(T_J^s,\der
T_J^s)$:

\begin{theorem}
Let $k$ be an algebraically closed field of characteristic $3$, let
$S$ be a para-Hurwitz algebra over $k$, let $J$ be the Jordan
algebra of $3\times 3$ hermitian matrices over the associated
Hurwitz algebra, and let $T_J^s$ be the symplectic triple system
attached to $J$ as in \eqref{eq:TJs}. Then the Lie superalgebras
$\frg(S,S_{4,2})$ and $\frtg(T_J^s,\der T_J^s)$ are isomorphic.
\end{theorem}

Therefore, the Lie superalgebras $\frg(S,S_{4,2})$ in the extended
Freudenthal Magic Square, for a para-Hurwitz algebra $S$, are
essentially the new simple Lie superalgebras in \cite[Theorem
3.2]{Eld3}. More specifically:

\begin{corollary}
Let $S_r$ ($r=1,2,4,8$) denote the unique para-Hurwitz algebra of
dimension $r$ over an algebraically closed field $k$ of
characteristic $3$. Then:
\begin{itemize}

\item $\frg(S_1,S_{4,2})$ is isomorphic to the simple Lie
superalgebra in \cite[Theorem 3.2(ii)]{Eld3}, obtained as
$\frg(T_J^s)$, for $J=H_3(k,*)$.

\item $[\frg(S_2,S_{4,2}),\frg(S_2,S_{4,2})]$ is isomorphic to the simple Lie
superalgebra in \cite[Theorem 3.2(iii)]{Eld3}, obtained as
$\frg(T_J^s)$, for $J=H_3(k\times k,*)$.

\item $\frg(S_4,S_{4,2})$ is isomorphic to the simple Lie
superalgebra in \cite[Theorem 3.2(iv)]{Eld3}, obtained as
$\frg(T_J^s)$, for $J=H_3(\Mat_2(k),*)$.

\item $\frg(S_8,S_{4,2})$ is isomorphic to the simple Lie
superalgebra in \cite[Theorem 3.2(v)]{Eld3}, obtained as
$\frg(T_J^s)$, for $J=H_3(C(k),*)$.

\end{itemize}
\end{corollary}


\section{Final remarks}

In the last sections, the Lie superalgebras $\frg(S_r,S_{1,2})$ and
$\frg(S_r,S_{4,2})$ ($r=1,2,4,8$) in the extended Freudenthal Magic
Square have been shown to be related to Lie superalgebras previously
constructed in terms of orthogonal and symplectic triple systems in
\cite{Eld3}.

Let us have a look in this last section to the remaining Lie
superalgebras in the extended square.

\smallskip

The result in \cite[Corollary 5.20]{CunEld1} shows that the even
part of $\frg(S_{1,2},S_{1,2})$ is isomorphic to the orthogonal Lie
algebra $\frso_7(k)$, while its odd part is the direct sum of two
copies of the spin module for $\frso_7(k)$, and therefore, over any
algebraically closed field of characteristic $3$,
$\frg(S_{1,2},S_{1,2})$ is isomorphic to the simple Lie superalgebra
in \cite[Theorem 4.23(ii)]{Eld3}, attached to a simple \emph{null
orthogonal triple system}.

Actually, the description of $\frg=\frg(S_{1,2},S_{1,2})$ in
\cite[\S 5.7]{CunEld1} shows that
\[
\begin{split}
\frg\subo&=\bigl(\frsp(V_1)\oplus\frsp(V_2)\oplus\frsl_2(k)\bigr)\oplus
\bigl(V_1\otimes V_2\otimes\frsl_2(k)\bigr),\\
\frg\subuno&=\bigl(V_1\oplus V_2\bigr)\otimes \frgl_2(k),
\end{split}
\]
and from here it is easy to directly check that $\frg\subo$ is
(isomorphic to) the orthogonal Lie algebra $\frso(\calV,\bil)$,
where $\calV=\bigl(V_1\otimes V_2\bigr)\oplus \frsl_2(k)$ and $\bil$
is the bilinear form (of maximal Witt index) on $\calV$ such that:
\[
\begin{split}
&\bil(V_1\otimes V_2,\frsl_2(k))=0,\\
&\bil(u_1\otimes u_2,v_1\otimes v_2)=\langle u_1\vert
v_1\rangle\langle u_2\vert v_2\rangle,\\
&\bil(p,q)=\det(p+q)-\det(p)-\det(q),
\end{split}
\]
for any $u_i,v_i\in V$ ($i=1,2$) and $p,q\in\frsl_2(k)$, and that
the Clifford algebra $Clif(\calV,-\bil)$ is isomorphic to
$\End_k(\frg\subuno)$, so that $\frg\subuno$ is isomorphic, as a
module for $\frg\subo\cong\frso(\calV,-\bil)\subseteq
Clif(\calV,-\bil)\subo$, to the direct sum of two copies of the spin
module (see \cite{Cunth} for details).

\smallskip

As for $\frg(S_{4,2},S_{4,2})$, the results in \cite[Proposition
5.10 and Corollary 5.11]{CunEld1} show that its even part is
isomorphic to the orthogonal Lie algebra $\frso_{13}(k)$, while its
odd part is the spin module for the even part, and hence that
$\frg(S_{4,2},S_{4,2})$ is the simple Lie superalgebra in
\cite[Theorem 3.1(ii)]{EldPac} for $l=6$.

Only the simple Lie superalgebra $\frg(S_{1,2},S_{4,2})$ has not
previously appeared in the literature. Its even part is isomorphic
to the symplectic Lie algebra $\frsp_8(k)$, while its odd part is
the irreducible module of dimension $40$ which appears as a
subquotient of the third exterior power of the natural module for
$\frsp_8(k)$ (see \cite[\S 5.5]{CunEld1}).

\smallskip

The results on orthogonal triple systems in Section 4, where
$\frg(S_{1,2},S)$ ($S$ a para-Hurwitz algebra) is shown to be
isomorphic to the Lie superalgebra $\frg(T)$ attached to the
orthogonal triple system defined on $J_0/k1$ for the Jordan algebra
of $3\times 3$ hermitian matrices on the Hurwitz algebra associated
to $S$, suggest the next definition in order to extend this
description of $\frg(S_{1,2},S)$ to the situation in which $S$ is a
para-Hurwitz superalgebra.

\begin{definition}
Let $T=T\subo\oplus T\subuno$ be a vector superspace over a field
$k$ of characteristic $\ne 2$ endowed with an even nonzero
supersymmetric bilinear form $(.\vert.):T\times T\rightarrow k$
(that is, $(T\subo\vert T\subuno)=0$, $(.\vert .)$ is symmetric on
$T\subo$ and alternating on $T\subuno$) and a triple product
$T\times T\times T\rightarrow T$: $(x,y,z)\mapsto [xyz]$
($[x_iy_jz_k]\in T_{i+j+k}$ for any $x_i\in T_i$, $y_j\in T_j$,
$z\in T_k$, $i,j,k=\bar 0$ or $\bar 1$). Then $T$ is said to be an
\emph{orthosymplectic triple system} if it satisfies the following
identities:
\begin{subequations}\label{eq:OSTS}
\begin{align}
&[xyz]+(-1)^{\lvert x\rvert\lvert y\rvert}[yxz]=0\label{eq:OSTSa}\\
&[xyz]+(-1)^{\lvert y\rvert\lvert z\rvert}[xzy]=
 (x\vert y)z+(-1)^{\lvert y\rvert\lvert z\rvert}(x\vert z)y
 -2(y\vert z)x\label{eq:OSTSb}\\
&[xy[uvw]]=[[xyu]vw]+
 (-1)^{(\lvert x\rvert+\lvert y\rvert)\lvert u\rvert}[u[xyv]w]\nonumber\\
 &\hspace{1.4in}+
 (-1)^{(\lvert x\rvert+\lvert y\rvert)(\lvert u\rvert
    +\lvert v\rvert)}[uv[xyw]]\label{eq:OSTSc}\\
&([xyu]\vert v)+(-1)^{(\lvert x\rvert+\lvert y\rvert)\lvert
         u\rvert}(u\vert [xyv])=0\label{eq:OSTSd}
\end{align}
\end{subequations}
for any homogeneous elements $x,y,u,v,w\in T$.
\end{definition}

The same arguments as in \cite[Theorem 2.9]{Eld3}, with suitable
parity signs here and there give:

\begin{theorem}\label{th:gTsortsy}
Let $T$ be an orthosymplectic triple system  and let
$\bigl(V,\langle.\vert.\rangle\bigr)$ be a two dimensional vector
space endowed with a nonzero alternating bilinear form. Let $\frs$
be a Lie subalgebra of $\der J$ containing $\inder T$. Define the
$\bZ_2$-graded superalgebra $\frg=\frg(T,\frs)=\frg(0)\oplus\frg(1)$
with
\[
\begin{cases}
\frg(0)= \frsp(V)\oplus \frs&\text{(so $\frg(0)\subo=\frsp(V)\oplus\frs\subo$,
$\frg(0)\subuno=\frs\subuno$),}\\
\frg(1)=V\otimes T&\text{(with $\frg(1)\subo=V\otimes T\subuno$,
$\frg(1)\subuno=V\otimes T\subo$, $V$ is odd!)},
\end{cases}
\]
and superanticommutative multiplication given by:
\begin{itemize}
\item
$\frg(0)$ is a subalgebra of $\frg$;
\item
$\frg(0)$ acts naturally on $\frguno$:
\[
[s,v\otimes x]=s(v)\otimes x,\qquad [d,v\otimes x]=(-1)^{\lvert
d\rvert}v\otimes d(x),
\]
for any $s\in \frsp(V)$, $d\in\frs$, $v\in V$, and $x\in T\subo\cup
T\subuno$;
\item
for any $u,v\in V$ and homogeneous $x,y\in T$:
\begin{equation}\label{eq:gTosproduct}
[u\otimes x,v\otimes y]=(-1)^{\lvert x\rvert}\bigl(
  -(x\vert y)\gamma_{u,v} +\langle u\vert v\rangle d_{x,y}\bigr)
\end{equation}
where $\gamma_{u,v}=\langle u\vert .\rangle v+\langle v\vert
.\rangle u$ and $d_{x,y}=[xy.]$.
\end{itemize}
Then $\frg(T,\frs)$ is a $\bZ_2$-graded Lie superalgebra. Moreover,
$\frg(T,\frs)$ is simple if and only if $\frs$ coincides with
$\inder T$ and $T$ is simple.

Conversely, given a $\bZ_2$-graded Lie superalgebra
$\frg=\frg(0)\oplus\frg(1)$ with
\[
\begin{cases}
\frg(0)=\frsp(V)\oplus \frs,\\
\frg(1)=V\otimes T,
\end{cases}
\]
where $T$ is a module for the superalgebra $\frs$ and $V$ is
considered as an odd vector space, by $\frsp(V)$-invariance of the
 bracket, equation \eqref{eq:gTosproduct} is satisfied for an even
supersymmetric bilinear form $(.\vert .):T\times T\rightarrow k$ and
a superskewsymmetric bilinear map $d_{.,.}:T\times T\rightarrow
\frs$. Then, if $(.\vert .)$ is not $0$ and a triple product on $T$
is defined by means of $[xyz]=d_{x,y}(z)$, $T$ becomes and
orthosymplectic triple system and the image of $\frs$ in $\frgl(T)$
under the given representation is a subalgebra of $\der T$
containing $\inder T$.
\end{theorem}

Now, given a para-Hurwitz superalgebra $S$, let $J$ be the Jordan
superalgebra of $3\times 3$ hermitian matrices over the associated
Hurwitz algebras as in Section 2. Consider, as in Section 4, the
quotient vector superspace
\begin{equation}\label{eq:TJos}
T_J^{os}=J_0/k1
\end{equation}
with even supersymmetric bilinear form induced by the trace:
\[
(\hat x\vert\hat y)=t(x\circ y)
\]
($\hat x=x+k1\in T_J^{os}$ for $x\in J_0$), and triple product given
by
\[
[\hat x\hat y\hat z]=\widehat{[L_x,L_y](z)}=\widehat{x\circ(y\circ
z)}-(-1)^{\lvert x\rvert\lvert y\rvert}\widehat{y\circ(x\circ z)}.
\]

Following the ideas in Section 4, consider the $\bZ_2$-graded
anticommutative superalgebra
\[
\frg=\frg(J)=\frg(0)\oplus\frg(1),
\]
with $\frg(0)=\frsp(V)\oplus\der J$, and $\frg(1)=V\otimes
T_J^{os}$.

The proof of Theorem \ref{th:gS12SgTJ}, with parity signs put all
over, works here to give that $\frg(S_{1,2},S)$ is isomorphic to
$\frg(J)$. Therefore, $\frg(J)$ is a Lie superalgebra!, because so
is $\frg(S_{1,2},S)$, and by Theorem \ref{th:gTsortsy}, $T_J^{os}$
is an orthosymplectic triple system, and $\frg(J)=\frg(T_J^{os},\der
J)$.

Therefore, there appears the following description of the Lie
superalgebras $\frg(S_{1,2},S)$ for any Hurwitz superalgebra $S$,
which includes a description of the Lie superalgebras
$\frg(S_{1,2},S_{1,2})$ and $\frg(S_{1,2},S_{4,2})$ in the extended
Freudenthal Magic Square \ref{ta:supersquare}.

\begin{theorem} Let $S$ be a para-Hurwitz superalgebra over a field
$k$ of characteristic $3$, and let $J$ be the Jordan superalgebra of
$3\times 3$ hermitian matrices over the associated Hurwitz
superalgebra. Then the Lie superalgebra $\frg(S_{1,2},S)$ is
isomorphic to the $\bZ_2$-graded Lie superalgebra
$\frg(J)=\frg(T_J^{os},\der J)$ associated to the simple
orthosymplectic triple system $T_J^{os}$ in \eqref{eq:TJos}.
\end{theorem}



\providecommand{\bysame}{\leavevmode\hbox
to3em{\hrulefill}\thinspace}
\providecommand{\MR}{\relax\ifhmode\unskip\space\fi MR }
\providecommand{\MRhref}[2]{%
  \href{http://www.ams.org/mathscinet-getitem?mr=#1}{#2}
} \providecommand{\href}[2]{#2}

\end{document}